\documentclass[reqno]{amsart}
\usepackage{mathrsfs}
\usepackage{color}
\usepackage{amsmath}
\usepackage{amsfonts}
\usepackage{amssymb}
\usepackage{graphicx}
\usepackage{hyperref}
\usepackage[numbers,sort&compress]{natbib}
\allowdisplaybreaks[4]

 \newtheorem{Theorem}{Theorem}[section]
 
 \newtheorem{Lemma}[Theorem]{Lemma}
 
 \newtheorem{Question}[Theorem]{Question}
 \newtheorem{Definition}[Theorem]{Definition}

 \newtheorem{Remark}[Theorem]{Remark}

 \numberwithin{equation}{section}

\begin{document}

\title[Optimal $L^2$ extension theorem on weakly pseudoconvex K\"{a}hler manifolds]
 {Guan-Zhou's unified version of optimal $L^2$ extension theorem on weakly pseudoconvex K\"{a}hler manifolds}

\author{Qi'an Guan}
\address{Qi'an Guan: School of
Mathematical Sciences, Peking University, Beijing 100871, China.}
\email{guanqian@math.pku.edu.cn}

\author{Zhitong Mi}
\address{Zhitong Mi: School of Mathematics and Statistics, Beijing Jiaotong University, Beijing
100044, China.
}
\email{zhitongmi@amss.ac.cn}
\author{Zheng Yuan}
\address{Zheng Yuan:  Institute of Mathematics, Academy of Mathematics and Systems Science, Chinese Academy of Sciences, Beijing 100190, China.}
\email{yuanzheng@amss.ac.cn}

\thanks{}

\subjclass[2020]{14F18, 32D15, 32U05, 32Q15}

\keywords{Holomorphic vector bundles, Optimal $L^2$ extension theorem, Weakly pseudoconvex manifolds}

\date{\today}

\dedicatory{}

\commby{}


\begin{abstract}
In this note, we establish Guan-Zhou's unified version of optimal $L^2$ extension theorem
 for holomorphic vector bundles with smooth hermitian metrics  on weakly pseudoconvex K\"{a}hler  manifolds. Combining with previous work of Guan-Mi-Yuan (\cite{GMYshm}), we generalize Guan-Zhou's unified version of optimal $L^2$ extension theorem to weakly pseudoconvex K\"{a}hler  manifolds.
\end{abstract}

\maketitle

\section{Introduction}

We recall the $L^2$ extension problem (see \cite{DemaillyAG}, see also \cite{zhou-abel}) as follows: let $Y$ be a complex subvariety of a complex manifold $M$; given a holomorphic object $f$ on $Y$ satisfying certain $L^2$ estimate on $Y$, finding a holomorphic extension $F$ of $f$ from $Y$ to $M$, together with a good or even optimal $L^2$ estimate of $F$ on $M$.

The existence part of $L^2$ extension problem was firstly solved by Ohsawa-Takegoshi \cite{OT87} and their result is called Ohsawa-Takegoshi $L^2$ extension theorem now.
After Ohsawa-Takegoshi $L^2$ extension theorem, many generalizations of $L^2$ extension theorem and applications of the theorem were established, e.g., see \cite{berndtsson annals,Demaillyshm, D2016,Ohsawa2,Ohsawa4,berndtsson paun,DHP,HPS}. Especially, $L^2$ extension theorems for different gain with non-optimal estimate were obtained by Berndtsson, Demailly, Ohsawa, Mcneal-Varolin, et al, see \cite{berndtsson1996,Ohsawa3,Ohsawa5,DemaillyManivel,McnealVarolin}.

The second part of $L^2$ extension problem was called the $L^2$ extension problem with optimal estimate or sharp $L^2$ extension problem (see \cite{zhou-abel}).
The method of undetermined functions was introduced to study the sharp $L^2$ extension problem by Guan-Zhou-Zhu \cite{GZZCRMATH,ZGZ}. For bounded
pseudoconvex domains in $\mathbb{C}^n$, Blocki \cite{Blocki-inv} developed the equation in the method of undetermined functions of Guan-Zhou-Zhu \cite{ZGZ}, and obtained the optimal version of Ohsawa-Takegoshi's $L^2$ extension theorem in \cite{OT87}. Using the method of undetermined functions,
Guan-Zhou \cite{guan-zhou13ap} (see also \cite{GZsci,guan-zhou CRMATH2012}) established an $L^2$ extension theorem with optimal estimate on Stein manifolds for continuous gain, which implied various optimal versions of $L^2$ extension theorem.  In \cite{ZhouZhu-jdg} and \cite{ZZ2019}, Zhou-Zhu proved optimal $L^2$ extension theorem  for  smooth gain on weakly pseudoconvex K\"{a}hler  manifolds.

Recall that, in  \cite{guan-zhou13ap} and  \cite{ZZ2019} (see also \cite{ZhouZhu-jdg}),  both the line bundle case (with singular metric) and the vector bundles case (with  smooth metric) were considered respectively. Note that \cite{ZZ2019} (see also \cite{ZhouZhu-jdg}) did not fully generalize the results in \cite{guan-zhou13ap}.
It is natural to ask the following question (we posed following question for  line bundles in \cite{GMYshm}).
\begin{Question}
\label{que1.1}Can one give unified versions of the optimal $L^2$ extension theorems (for both line bundles with singular metrics and vector bundles with  smooth metrics) of Guan-Zhou \cite{guan-zhou13ap} and Zhou-Zhu \cite{ZZ2019} (see also \cite{ZhouZhu-jdg}).
\end{Question}

 In \cite{GMYshm}, we presented an optimal $L^2$ extension theorems for holomorphic vector bundles equipped with singular Nakano-semipositive metrics on weakly pseudoconvex K\"{a}hler  manifolds, which can be viewed as a unified version of  the line bundle case (with singular metric) in \cite{guan-zhou13ap} and \cite{ZZ2019} and hence answered Question \ref{que1.1} for line bundles (with singular metrics).

In this note, we give an affirmative answer to Question \ref{que1.1} for the  vector bundle  case  (with  smooth metric) by presenting an optimal $L^2$ extension theorem of holomorphic vector bundles  with smooth hermitian metrics for continuous gain on weakly pseudoconvex K\"{a}hler  manifolds, which is a unified version of the optimal $L^2$ extension theorems for vector bundles (with  smooth metrics) of \cite{guan-zhou13ap} and \cite{ZZ2019}.

\subsection {Main result}

\begin{Definition}
	\label{def:neat}A function $\psi:M\rightarrow[-\infty,+\infty)$ on a complex manifold $M$ is said to be quasi-plurisubharmonic if $\psi$ is locally the sum of a plurisubharmonic function and a smooth function (or equivalently, if $i\partial\bar\partial\psi$ is locally bounded from below). In addition, we say that $\psi$ has neat analytic singularities if every point $z\in M$ possesses an open neighborhood $U$ on which $\psi$ can be written as
	$$\psi=c\log\sum_{1\le j\le N}|g_j|^2+v,$$
	where $c\ge0$ is a constant, $g_j\in\mathcal{O}(U)$ and $v\in C^{\infty}(U)$.
\end{Definition}

\begin{Definition}
	If $\psi$ is a quasi-plurisubharmonic function on an $n$-dimensional complex manifold $M$, the multiplier ideal sheaf $\mathcal{I}(\psi)$ is the coherent analytic subsheaf of $\mathcal{O}_M$ defined by
	\begin{displaymath}
		\mathcal{I}(\psi)_z=\left\{f\in\mathcal{O}_{M,z}:\exists U\ni z,\,\int_U|f|^2e^{-\psi}d\lambda<+\infty \right\},
	\end{displaymath}
			where $U$ is an open coordinate neighborhood of $z$ and $d\lambda$ is the Lebesgue measure in the corresponding open chart of $\mathbb{C}^n$.
			
			We say that the singularities of $\psi$ are log canonical along the zero variety $Y=V(I(\psi))$ if $\mathcal{I}((1-\epsilon)\psi)|_Y=\mathcal{O}_M|_Y$ for any $\epsilon>0$.
\end{Definition}

Let $(M,\omega)$ be an $n$-dimensional K\"ahler manifold, and let $dV_{M,\omega}=\frac{1}{n!}\omega^n$ be the corresponding K\"ahler volume element.

\begin{Definition}
	\label{def:oshawa measure}Let $\psi$ be a quasi-plurisubharmonic function on $M$ with neat analytic singularities. Assume that the singularities of $\psi$ are log canonical along the zero variety $Y=V(I(\psi))$. Denote $Y^0=Y_{\rm{reg}}$ the regular point set of $Y$. If $g\in C_c(Y^0)$  and $\hat g\in C_c(M)$ satisfy $\hat g|_{Y^0}=g$ and $(supp\,\hat{g})\cap Y=Y^0$, we set
	\begin{equation}
\label{eq:0627e}
		\int_{Y^0}gdV_{M,\omega}[\psi]=\limsup_{t\rightarrow+\infty}\int_{\{-t-1<\psi<-t\}}\hat ge^{-\psi}dV_{M,\omega}.
	\end{equation}
\end{Definition}
\begin{Remark}[see \cite{D2016}]\label{r:measure}
	By Hironaka's desingularization theorem, it is not hard to see that the limit in the right of equality \eqref{eq:0627e} does not depend on the continuous extension $\hat g$ and $dV_{M,\omega}[\psi]$ is well defined on $Y^0$ .
\end{Remark}

We recall the following class of functions, so called ``gain''.

\begin{Definition}[see \cite{guan-zhou13ap}]
\label{class gTdelta}Let $T\in (-\infty,+\infty)$ and $\delta\in (0,+\infty)$. Let $\mathcal{G}_{T,\delta}$ be the class of functions $c(t)$ which satisfies the following statements,\\
(1) $c(t)$ is a continuous positive function on $[T,+\infty)$,\\
(2)  $\int_{T}^{+\infty}c(t)e^{-t}dt<+\infty$,\\
(3) for any $t> T$, the following equality holds,
\begin{equation}\label{integarl condition 1}
\begin{split}
    &{\left(\frac{1}{\delta}c(T)e^{-T}+\int_{T}^{t}c(t_1)e^{-t_1}dt_1\right)}^2 > \\
     &c(t)e^{-t}\bigg(\int_{T}^{t}(\frac{1}{\delta}c(T)e^{-T}+\int_{T}^{t_2}c(t_1)e^{-t_1}dt_1)dt_2
     +\frac{1}{{\delta}^2}c(T)e^{-T}\bigg).
\end{split}
\end{equation}
\end{Definition}
\begin{Remark}The number $-T$, $\frac{1}{\delta}$ and function $c(t)$ are equal to the number $\alpha_0$, $\alpha_1$ and function $\frac{1}{R(-t)e^{-t}}$ in \cite{ZZ2019}. We would like to use $-T$, $\frac{1}{\delta}$ and $c(t)$ in this note.
\end{Remark}

In this note, we establish Guan-Zhou's unified version of optimal $L^2$ extension theorem
 for holomorphic vector bundles with smooth hermitian metrics  on weakly pseudoconvex K\"{a}hler  manifolds.

\begin{Theorem} [Main theorem] \label{main result}
Let $c(t)\in \mathcal{G}_{T,\delta}$, where $\delta<+\infty$. Let $(M,\omega)$ be a weakly pseudoconvex K\"ahler manifold. Let $\psi<-T$ be a quasi-plurisubharmonic function on $M$ with neat analytic singularities. Let $Y:=V(\mathcal{I}(\psi))$ and assume that $\psi$ has log canonical singularities along $Y$.  Let $E$ be a rank $r$ holomorphic vector bundle over $M$   equipped with a smooth Hermitian metric $h$.  Assume that \\
(1) $\sqrt{-1}\Theta_{h}+\sqrt{-1}\partial\bar{\partial}\psi\ge 0$ on $M\backslash\{\psi=-\infty\}$ in the sense of Nakano;\\
(2)  $\sqrt{-1}\Theta_{h}+\sqrt{-1}\partial\bar{\partial}\psi
+\frac{1}{s(-\psi)}\sqrt{-1}\partial\bar{\partial}\psi\ge 0$ on $M\backslash\{\psi=-\infty\}$ in the sense of Nakano, where
$$s(t):=\frac{\int^t_T\bigg(\frac{1}{\delta}c(T)e^{-T}+\int^{t_2}_T c(t_1)e^{-t_1}dt_1\bigg)dt_2+\frac{1}{\delta^2}c(T)e^{-T}}{\frac{1}{\delta}c(T)e^{-T}+\int^t_T
c(t_1)e^{-t_1}dt_1}.$$

Then for every section $f \in H^0(Y^0,(K_M\otimes E)|_{Y^0})$ on $Y^0=Y_{\text{reg}}$ such that
\begin{equation}\label{mainth:ohsawa measure finite}
  \int_{Y^0}|f|^2_{\omega,h}dV_{M,\omega}[\psi]<+\infty,
\end{equation}
there exists a section $F\in H^0(M,K_M\otimes E)$ such that $F|_{Y_0}=f$ and
\begin{equation}\label{mainth:l2 estimate}
 \int_M c(-\psi)|F|^2_{\omega,h}dV_{M,\omega}\le \left(\frac{1}{\delta}c(T)e^{-T}+\int_{T}^{+\infty}c(t_1)e^{-t_1}dt_1\right)\int_{Y^0}|f|^2_{\omega,h}dV_{M,\omega}[\psi].
\end{equation}
\end{Theorem}
Combining Theorem \ref{main result} with previous work of Guan-Mi-Yuan (see \cite{GMYshm}), we generalize Guan-Zhou's unified version of optimal $L^2$ extension theorem to weakly pseudoconvex K\"{a}hler  manifolds.
\begin{Remark}
When $M$ is a Stein manifold, Theorem \ref{main result} can be referred to Guan-Zhou \cite{guan-zhou13ap}. When $c(t)$ is a smooth function, $\liminf_{t\to +\infty}c(t)>0$ and $c(t)e^{-t}$ is decreasing with respect to $t$ near $+\infty$, Theorem \ref{main result} can be referred to Zhou-Zhu \cite{ZZ2019} (see also \cite{ZhouZhu-jdg}).
\end{Remark}

\section{Preparations}

\subsection{Preparations for the proof of main theorem}
\
In this section, we make some preparations for the proof of main theorem.

We would like to recall some lemmas which will be used in this section.
\begin{Lemma}
[Theorem 1.5 in \cite{Demailly82}]
\label{completeness}
Let M be a K\"ahler manifold, and Z be an analytic subset of M. Assume that
$\Omega$ is a relatively compact open subset of M possessing a complete K\"ahler
metric. Then $\Omega\backslash Z $ carries a complete K\"ahler metric.

\end{Lemma}

\begin{Lemma}
[Lemma 6.9 in \cite{Demailly82}]
\label{extension of equality}
Let $\Omega$ be an open subset of $\mathbb{C}^n$ and Z be a complex analytic subset of
$\Omega$. Assume that $v$ is a (p,q-1)-form with $L^2_{loc}$ coefficients and h is
a (p,q)-form with $L^1_{loc}$ coefficients such that $\bar{\partial}v=h$ on
$\Omega\backslash Z$ (in the sense of distribution theory). Then
$\bar{\partial}v=h$ on $\Omega$.
\end{Lemma}

\begin{Lemma}
[Remark 3.2 in \cite{DemaillyManivel}]
\label{d-bar equation with error term}
Let ($M,\omega$) be a complete K\"ahler manifold equipped with a (non-necessarily
complete) K\"ahler metric $\omega$, and let Q be a Hermitian vector bundle over M.
Assume that $\eta$ and g are smooth bounded positive functions on M and let
$B:=[\eta \sqrt{-1}\Theta_Q-\sqrt{-1}\partial \bar{\partial} \eta-\sqrt{-1}g
\partial\eta \wedge\bar{\partial}\eta, \Lambda_{\omega}]$. Assume that $\delta \ge
0$ is a number such that $B+\delta I$ is semi-positive definite everywhere on
$\wedge^{n,q}T^*M \otimes Q$ for some $q \ge 1$. Then given a form $v \in
L^2(M,\wedge^{n,q}T^*M \otimes Q)$ such that $D^{''}v=0$ and $\int_M \langle
(B+\delta I)^{-1}v,v\rangle_Q dV_M < +\infty$, there exists an approximate solution
$u \in L^2(M,\wedge^{n,q-1}T^*M \otimes Q)$ and a correcting term $h\in
L^2(M,\wedge^{n,q}T^*M \otimes Q)$ such that $D^{''}u+\sqrt{\delta}h=v$ and
\begin{equation}
\int_M(\eta+g^{-1})^{-1}|u|^2_QdV_M+\int_M|h|^2_QdV_M \le \int_M \langle (B+\delta
I)^{-1}v,v\rangle_Q dV_M.
\end{equation}
\end{Lemma}
Let $M$ be a complex manifold. Let $\omega$ be a continuous hermitian metric on $M$. Let $dV_M$ be a continuous volume form on $M$. We denote by $L^2_{p,q}(M,\omega,dV_M)$ the spaces of $L^2$ integrable $(p,q)$-forms over $M$ with respect to $\omega$ and $dV_M$. It is known that $L^2_{p,q}(M,\omega,dV_M)$ is a Hilbert space.
\begin{Lemma}[see Lemma 9.1 in \cite{GMY-boundary5}]
\label{weakly convergence}
Let $\{u_n\}_{n=1}^{+\infty}$ be a sequence of $(p,q)$-forms in $L^2_{p,q}(M,\omega,dV_M)$ which is weakly convergent to $u$. Let $\{v_n\}_{n=1}^{+\infty}$ be a sequence of Lebesgue measurable real functions on $M$ which converges pointwisely to $v$. We assume that there exists a constant $C>0$ such that $|v_n|\le C$ for any $n$. Then $\{v_nu_n\}_{n=1}^{+\infty}$ weakly converges to $vu$ in $L^2_{p,q}(M,\omega,dV_M)$.
\end{Lemma}

Let $X$ be an $n-$dimensional complex manifold and $\omega$ be a hermitian metric on $X$. Let $Q$ be a holomorphic vector bundle on $X$ with rank $r$. Let $\{h_i\}_{i=1}^{+\infty}$ be a family of $C^2$ smooth hermitian metric on $Q$ and $h$ be a measurable metric on $Q$ such that $\lim_{i\to+\infty}h_i=h$ almost everywhere on $X$.  We assume that
$h_i$ is increasingly convergent to $h$ as $i\to+\infty$.

The following optimal $L^2$ extension theorem for vector bundles with smooth hermitian metric on Stein manifolds will be used in our discussion.
\begin{Theorem}[see \cite{guan-zhou13ap}]\label{l2 extension on stein} Let $c(t)\in \mathcal{G}_{T,\delta}$ for some $T\in(-\infty,+\infty)$ and $0<\delta<+\infty$.
Let $M$ be a Stein manifold and $\omega$ be a hermitian metric on $M$. Let $h$ be a smooth hermitian metric on a holomorphic vector bundle $E$ on $M$ with rank $r$. Let $\psi<-T$ be a quasi-plurisubharmonic function on $X$ with neat analytic singularities. Let $Y:=V(\mathcal{I}(\psi))$ and assume that $\psi$ has log canonical singularities along $Y$. Assume that\\
(1) $\sqrt{-1}\Theta_h+\sqrt{1}\partial\bar{\partial}\psi$ is Nakano semi-positive on $M\backslash\{\psi=-\infty\}$, \\
(2) there exists a continuous function $a(t)$ on $(T,+\infty]$ such that $0<a(t)\le s(t)$ and $a(-\psi)\sqrt{-1}\Theta_{he^{-\psi}}+\sqrt{1}\partial\bar{\partial}\psi$ is Nakano semi-positive on $M\backslash\{\psi=-\infty\}$, where
$$s(t):=\frac{
  \int_{T}^{t}\big(\frac{1}{\delta}c(T)e^{-T}+\int_{T}^{t}c(t_1)e^{-t_1}dt_1\big)dt_2+\frac{1}{\delta^2}c(T)e^{-T}}
  {\frac{1}{\delta}c(T)e^{-T}+\int_{T}^{t}c(t_1)e^{-t_1}dt_1}.$$
Then  for any holomorphic section $f$ of $K_M\otimes E|_Y$ on $Y$ satisfying
$$\int_{Y_0}|f|^2_{\omega,h}dV_{M,\omega}[\psi]<+\infty,$$
there exists a holomorphic section $F$ of   $K_M\otimes E$ on $M$ satisfying $F|_Y=f$ and
$$\int_M c(-\psi)|F|^2_{\omega,h}dV_{M,\omega}\le \left(\frac{1}{\delta}c(T)e^{-T}+\int_{T}^{+\infty}c(t_1)e^{-t_1}dt_1\right)\int_{Y_0}|f|^2_{\omega,h}dV_{M,\omega}[\psi].$$

\end{Theorem}

The following lemma will be used in the proof of the main theorem.
\begin{Lemma}[see Theorem 4.4.2 in \cite{Hormander}] \label{hormander} Let $\Omega$ be a pseudoconvex domain in $\mathbb{C}^n$, and $\varphi$ be a plurisubharmonic function on $\Omega$. For any $w\in L^2_{p,q+1}(\Omega,e^{-\varphi})$ with $\bar{\partial}w=0$, there exists a solution $s\in L^2_{p,q}(\Omega,e^{-\varphi})$ of the equation $\bar{\partial}s=w$ such that
$$\int_{\Omega}\frac{|s|^2}{(1+|z|^2)^2}e^{-\varphi}d\lambda\le \int_{\Omega}|w|^2e^{-\varphi}d\lambda,$$
where $d\lambda$ is the Lebesgue measure on $\mathbb{C}^n$.
\end{Lemma}

Let $c(t)$  belong to class $\mathcal{G}_{T,\delta}$.
Recall that
$$s(t):=\frac{\int^t_T\bigg(\frac{1}{\delta}c(T)e^{-T}+\int^{t_2}_T c(t_1)e^{-t_1}dt_1\bigg)dt_2+\frac{1}{\delta^2}c(T)e^{-T}}{\frac{1}{\delta}c(T)e^{-T}+\int^t_T
c(t_1)e^{-t_1}dt_1}.$$
We have following regularization lemma for $c(t)$.
\begin{Lemma}[see \cite{GMYshm}]
\label{approx of ct}
Let $c(t)\in \mathcal{G}_{T,\delta}$. Let $\{\beta_m<1\}$ be a sequence of positive real numbers such that $\beta_m$ decreasingly converges to $0$ as $m\to +\infty$. Then there exists a sequence of positive functions $c_m(t)$ on $[T,+\infty)$, which satisfies:\\
(1) $c_m(t)\in \mathcal{G}_{T,\delta}$;\\
(2) $c_m(t)e^{-t}$ is decreasing with respect to $t$ near $+\infty$;\\
(3) $c_m(t)$ is smooth on $[T+4\beta_m,+\infty)$;\\
(4) $c_m(t)$ are uniformly convergent to $c(t)$ on any compact subset of $(T,+\infty)$;\\
(5) $\frac{1}{\delta}c_m(T)e^{-T}+\int_{T}^{+\infty}c_m(t)e^{-t}dt$ converges to $\frac{1}{\delta}c(T)e^{-T}+\int_{T}^{+\infty}c(t)e^{-t}dt<+\infty$ as $m\to +\infty$;\\
(6) For each $m$, there exists $\kappa_m>0$ such that
$$S_m(t):=\frac{\int^t_T\bigg(\frac{1}{\delta}c_m(T)e^{-T}+\int^{t_2}_T c_m(t_1)e^{-t_1}dt_1\bigg)dt_2+\frac{1}{\delta^2}c_m(T)e^{-T}+\kappa_m}{\frac{1}{\delta}c_m(T)e^{-T}+\int^t_T
c_m(t_1)e^{-t_1}dt_1}>s(t),$$
for any $t\ge T$ and $S'_m(t)>0$ on $[T+\beta_m,+\infty)$.
\end{Lemma}

The following Lemma will be used in the proof of the Theorem \ref{main result}.
\begin{Lemma}[see \cite{GY-concavity1}]
	\label{l:converge}
	Let $M$ be a complex manifold. Let $S$ be an analytic subset of $M$.  	
	Let $\{g_j\}_{j=1,2,...}$ be a sequence of nonnegative Lebesgue measurable functions on $M$, which satisfies that $g_j$ are almost everywhere convergent to $g$ on  $M$ when $j\rightarrow+\infty$,  where $g$ is a nonnegative Lebesgue measurable function on $M$. Assume that for any compact subset $K$ of $M\backslash S$, there exist $s_K\in(0,+\infty)$ and $C_K\in(0,+\infty)$ such that
	$$\int_{K}{g_j}^{-s_K}dV_M\leq C_K$$
	 for any $j$, where $dV_M$ is a continuous volume form on $M$.
	
 Let $\{F_j\}_{j=1,2,...}$ be a sequence of holomorphic $(n,0)$ form on $M$. Assume that $\liminf_{j\rightarrow+\infty}\int_{M}|F_j|^2g_j\leq C$, where $C$ is a positive constant. Then there exists a subsequence $\{F_{j_l}\}_{l=1,2,...}$, which satisfies that $\{F_{j_l}\}$ is uniformly convergent to a holomorphic $(n,0)$ form $F$ on $M$ on any compact subset of $M$ when $l\rightarrow+\infty$, such that
 $$\int_{M}|F|^2g\leq C.$$
\end{Lemma}

\section{Proof of Theorem \ref{main result}}
In this section, we prove Theorem \ref{main result}.

\begin{proof}

  As $M$ is weakly pseudoconvex, there exists a smooth plurisubharmonic
exhaustion function $P$ on $M$. Let $M_k:=\{P<k\}$ $(k=1,2,...,) $. We choose $P$ such that
$M_1\ne \emptyset$.\par
Then $M_k$ satisfies $
M_k\Subset  M_{k+1}\Subset  ...M$ and $\cup_{k=1}^n M_k=M$. Each $M_k$ is weakly
pseudoconvex K\"ahler manifold with exhaustion plurisubharmonic function
$P_k=1/(k-P)$.
\par
\emph{We will fix $k$ during our discussion until the last step.}

\

\textbf{Step 1: regularization of $c(t)$.}

\

As $e^{\psi}$ is a smooth function on $M$ and $\psi<-T$ on $M$, we know that
$$\sup_{M_k}\psi<-T-8\epsilon_k,$$
where $\epsilon_k >0$ is a real number depending on $k$.

It follows from $c(t)$ belongs to class $\mathcal{G}_{T,\delta}$, by Lemma \ref{approx of ct}, that we have a sequence of functions $\{c_{k}(t)\}_{k\in\mathbb{Z}^+}$ which satisfies $c_k(t)$ is continuous on $[T,+\infty)$, and smooth on $[T+4\epsilon_k,+\infty)$ and other conditions in Lemma \ref{approx of ct}. Condition (6) of Lemma \ref{approx of ct} tells that
$$S_k(t):=\frac{\int^t_T\bigg(\frac{1}{\delta}c_k(T)e^{-T}+\int^{t_2}_T c_k(t_1)e^{-t_1}dt_1\bigg)dt_2+\frac{1}{\delta^2}c_k(T)e^{-T}+\kappa_k}{\frac{1}{\delta}c_k(T)e^{-T}+\int^t_T
c_k(t_1)e^{-t_1}dt_1}>s(t),$$
for any $t\ge T$ and $S'_k(t)>0$ on $[T+\epsilon_k,+\infty)$.

 As $S_k(t)>S(t)$ on  $t\ge T$, we know that
 $$S_k(-\psi)\big(\sqrt{-1}\partial\bar{\partial}\varphi+\sqrt{-1}\partial\bar{\partial}\psi\big)
+\sqrt{-1}\partial\bar{\partial}\psi\ge 0$$ on $M\backslash\{\psi=-\infty\}$ in the sense of currents. Denote $u_k(t):=-\log(\frac{1}{\delta}c(T)e^{-T}+\int_{T}^{t}c_k(t_1)e^{-t_1}dt_1)$. We note that we still have $S'_k(t)-S_k(t)u'_k(t)=1$ and $(S_k(t)+\frac{S'^2_k(t)}{u''_k(t)S_k(t)-S''_k(t)})e^{u_k(t)-t}=\frac{1}{c_k(t)}$.

\

\textbf{Step 2: construction of a family of smooth extensions $\tilde{f}$ of $f$ to a neighborhood of $ Y$ to $M$.}

\

Let $\mathcal{U}=\{U_i\}_{i\in I}$ be a locally finite covering  of $M$ by local coordinate balls.
It follows from Theorem \ref{l2 extension on stein} that there exists a holomorphic section $f_i\in \Gamma\big(U_i,\mathcal{O}_M(K_M\otimes E)\big)$ which is a holomorphic extension of $f$ from $U_i\cap Y^0$ to $U_i$. Let $\{\xi_i\}_{i\in I}$ be a partition of unity subordinate to $\mathcal{U}$, and denote
$$\tilde{f}:=\sum_{i\in I}\xi_i f_i.$$
Then $\tilde{f}$ is smooth on $M$, and we have
\begin{equation}\label{d-baroftildef}
  \begin{split}
     \bar{\partial}\tilde{f}|_{U_j}& =\bar{\partial}\tilde{f}-\bar{\partial}f_j \\
       & =\bar{\partial}(\sum_{i\in I}\xi_i f_i)-\bar{\partial}(\sum_{i\in I}\xi_i f_j)\\
       &=\sum_{i\in I}\bar{\partial} \xi_i\wedge (f_i-f_j), \text{ for all }j\in I.
  \end{split}
\end{equation}

Note that $f_i-f_j=0$ on $U_i\cap U_j\cap Y$. It follows from $\psi$ has neat analytic singularities and  is log canonical along the zero variety $Y$ that we know that
\begin{equation}\label{fi-fj}
f_i-f_j\in \Gamma\big(U_i\cap U_j,\mathcal{O}_M(K_M\otimes E)\otimes \mathcal{I}(\psi)\big).
\end{equation}

Note that $k$ is fixed until the last step and $M_k\Subset M$. We may assume that $M_k\subset \cup_{i=1}^{N}U_i$, where $N$ is a positive integer.

\

\textbf{Step 3: recall some notations.}

\

Let $\epsilon \in (0,\frac{1}{8})$. Let $\{v_{t_0,\epsilon}\}_{\epsilon \in
(0,\frac{1}{8})}$ be a family of smooth increasing convex functions on $\mathbb{R}$, such
that:
\par
(1) $v_{t_0,\epsilon}(t)=t$ for $t\ge-t_0-\epsilon$, $v_{\epsilon}(t)=constant$ for
$t<-t_0-1+\epsilon$;\par
(2) $v_{\epsilon}{''}(t)$ are convergence pointwisely
to $\mathbb{I}_{(-t_0-1,-t_0)}$,when $\epsilon \to 0$, and $0\le
v_{\epsilon}{''}(t) \le \frac{1}{1-4\epsilon}\mathbb{I}_{(-t_0-1+\epsilon,-t_0-\epsilon)}$
for ant $t \in \mathbb{R}$;\par
(3) $v_{\epsilon}{'}(t)$ are convergence pointwisely to $b(t)$ which is a continuous
function on $\mathbb{R}$ when $\epsilon \to 0$ and $0 \le v_{\epsilon}{'}(t) \le 1$ for any
$t\in \mathbb{R}$.\par
One can construct the family $\{v_{t_0,\epsilon}\}_{\epsilon \in (0,\frac{1}{8})}$  by
 setting
\begin{equation}\nonumber
\begin{split}
v_{t_0,\epsilon}(t):=&\int_{-\infty}^{t}\bigg(\int_{-\infty}^{t_1}(\frac{1}{1-4\epsilon}
\mathbb{I}_{(-t_0-1+2\epsilon,-t_0-2\epsilon)}*\rho_{\frac{1}{4}\epsilon})(s)ds\bigg)dt_1\\
&-\int_{-\infty}^{-t_0}\bigg(\int_{-\infty}^{t_1}(\frac{1}{1-4\epsilon}
\mathbb{I}_{(-t_0-1+2\epsilon,-t_0-2\epsilon)}*\rho_{\frac{1}{4}\epsilon})(s)ds\bigg)dt_1-t_0,
\end{split}
\end{equation}
where $\rho_{\frac{1}{4}\epsilon}$ is the kernel of convolution satisfying
$\text{supp}(\rho_{\frac{1}{4}\epsilon})\subset
(-\frac{1}{4}\epsilon,{\frac{1}{4}\epsilon})$.
Then it follows that
\begin{equation}\nonumber
v_{t_0,\epsilon}{''}(t)=\frac{1}{1-4\epsilon}
\mathbb{I}_{(-t_0-1+2\epsilon,-t_0-2\epsilon)}*\rho_{\frac{1}{4}\epsilon}(t),
\end{equation}
and
\begin{equation}\nonumber
v_{t_0,\epsilon}{'}(t)=\int_{-\infty}^{t}\bigg(\frac{1}{1-4\epsilon}
\mathbb{I}_{(-t_0-1+2\epsilon,-t_0-2\epsilon)}*\rho_{\frac{1}{4}\epsilon}\bigg)(s)ds.
\end{equation}
Note that $\text{supp}v_{t_0,\epsilon}{''}(t)\Subset (-t_0-1+\epsilon,-t_0-\epsilon)$ and
$\text{supp}\big(1- v_{t_0,\epsilon}{'}(t)\big)\Subset (-\infty,-t_0-\epsilon)$

We also note that $S_k\in C^{\infty}([T+4\epsilon_k,+\infty))$ satisfies
 $S_k'>0$ on $[T+\epsilon_k,+\infty)$ and $u_k\in C^{\infty}([T+4\epsilon_k,+\infty))$ satisfies
$\lim_{t\to +\infty}u_k(t)=-\log(\frac{1}{\delta}c_k(T)e^{-T}+\int^{+\infty}_T c_k(t_1)e^{-t_1}dt_1)$ and $u'_k<0$.
Recall that $u_k(t)$ and $S_k(t)$ satisfy $$S'_k(t)-S_k(t)u'_k(t)=1$$ and $$(S_k(t)+\frac{S'^2_k(t)}{u''_k(t)S_k(t)-S''_k(t)})e^{u_k(t)-t}=\frac{1}{c_k(t)}.$$
Note that $u_k''S_k-S_k''=-S'_ku'_k>0$ on $[T+2\epsilon_k,+\infty)$.
Denote $\tilde{g}_k(t):=\frac{u''_kS_k-S''_k}{S'^2_k}(t)$, then $\tilde{g}_k(t)$ is a positive smooth function on $[T+4\epsilon_k,+\infty)$.

Denote $\eta:=S_k(-v_{t_0,\epsilon}(\psi))$, $\phi:=u_k(-v_{t_0,\epsilon}(\psi))$ and $g:=\tilde{g}_k(-v_{t_0,\epsilon}(\psi))$. Then $\eta$ and $g$ are smooth bounded positive functions on $M_k$ such that $\eta+g^{-1}$ is a smooth bounded positive function on $M_k$.

Denote $\sum:=\{\psi=-\infty\}$.
As $\psi$ has neat analytic singularities, we know that $\sum$ is an analytic subset of $M$ and $\psi$ is smooth on $M\backslash \sum$. Note that, by Lemma \ref{completeness}, $M_k\backslash \sum$ carries a complete K\"ahler metric.
 Denote $\tilde{h}:=he^{-\psi}e^{-\phi}$ on $M_k\backslash  \sum$.

\

\textbf{Step 4: some calculations.}\label{step 4 of main result}

\

We set
$B=[\eta \sqrt{-1}\Theta_{\tilde{h}}-\sqrt{-1}\partial \bar{\partial}
\eta\otimes\text{Id}_E-\sqrt{-1}g\partial\eta \wedge\bar{\partial}\eta\otimes\text{Id}_E, \Lambda_{\omega}]$ on $M_k\backslash \sum$. Direct calculation shows that
\begin{equation}\nonumber
\begin{split}
\partial\bar{\partial}\eta=&
-S_k'\big(-v_{t_0,\epsilon}(\psi)\big)\partial\bar{\partial}\big(v_{t_0,\epsilon}(\psi)\big)
+S_k''\big(-v_{t_0,\epsilon}(\psi)\big)\partial\big(v_{t_0,\epsilon}(\psi)\big)\wedge
\bar{\partial}\big(v_{t_0,\epsilon}(\psi)\big),\\
\eta\Theta_{\tilde{h}}=&\eta\partial\bar{\partial}\phi\otimes\text{Id}_E+\eta\Theta_{h}+
\eta\partial\bar{\partial}\psi\otimes\text{Id}_E\\
=&S_ku_k''\big(-v_{t_0,\epsilon}(\psi)\big)\partial\big(v_{t_0,\epsilon}(\psi)\big)\wedge
\bar{\partial}\big(v_{t_0,\epsilon}(\psi)\big)\otimes\text{Id}_E
-S_ku_k'\big(-v_{t_0,\epsilon}(\psi)\big)\partial\bar{\partial}\big(v_{t_0,\epsilon}(\psi)\big)\otimes\text{Id}_E\\
+&S_k\Theta_{h}
+S_k\partial\bar{\partial}\psi\otimes\text{Id}_E.
\end{split}
\end{equation}
Hence
\begin{equation}\nonumber
\begin{split}
&\eta \sqrt{-1}\Theta_{\tilde{h}}-\sqrt{-1}\partial \bar{\partial}
\eta\otimes\text{Id}_E-\sqrt{-1}g\partial\eta \wedge\bar{\partial}\eta\otimes\text{Id}_E\\
=&S_k\sqrt{-1}\Theta_{h}
+S_k\sqrt{-1}\partial\bar{\partial}\psi\otimes\text{Id}_E\\
+&\big(S_k'-S_ku_k'\big)\big(v'_{t_0,\epsilon}(\psi)\sqrt{-1}\partial\bar{\partial}(\psi)+
v''_{t_0,\epsilon}(\psi)\sqrt{-1}\partial(\psi)\wedge\bar{\partial}(\psi)\big)\otimes\text{Id}_E\\
+&[\big(u_k''S_k-S_k''\big)-\tilde{g}_kS_k'^2]\sqrt{-1}\partial\big(v_{t_0,\epsilon}(\psi)\big)\wedge\bar{\partial}\big(v_{t_0,\epsilon}(\psi)\big)\otimes\text{Id}_E,
\end{split}
\end{equation}
where we omit the term $-v_{t_0,\epsilon}(\psi)$ in $(S_k'-S_ku_k')\big(-v_{t_0,\epsilon}(\psi)\big)$ and $[(u_k''S_k-S_k'')-\tilde{g}_kS_k'^2]\big(-v_{t_0,\epsilon}(\psi)\big)$ for simplicity.
Note that $S_k'(t)-S_k(t)u_k'(t)=1$, $\frac{u_k''(t)S_k(t)-S_k''(t)}{S_k'^2(t)}-\tilde{g}_k(t)=0$. We have
\begin{equation}\label{cal of curvature 0}
\begin{split}
&\eta \sqrt{-1}\Theta_{\tilde{h}}-\sqrt{-1}\partial \bar{\partial}
\eta\otimes\text{Id}_E-\sqrt{-1}g\partial\eta \wedge\bar{\partial}\eta\otimes\text{Id}_E\\
=&S_k\sqrt{-1}\Theta_{h}
+S_k\sqrt{-1}\partial\bar{\partial}\psi\otimes\text{Id}_E\\
+&\big(v'_{t_0,\epsilon}(\psi)\sqrt{-1}\partial\bar{\partial}(\psi)+
v''_{t_0,\epsilon}(\psi)\sqrt{-1}\partial(\psi)\wedge\bar{\partial}(\psi)\big)\otimes\text{Id}_E.
\end{split}
\end{equation}
We would like to discuss a property of $S_k(t)$.
\begin{Lemma}[see \cite{GMYshm}]
\label{lem:main}
For large enough $t_{0}$, and for any $\varepsilon\in(0,1/4)$, the inequality
\begin{equation}
\label{equ:main}
S_k(-v_{t_{0},\varepsilon}(t))\geq S_k(-t)v'_{t_{0},\varepsilon}(t)
\end{equation}
holds for any $t\in(-\infty,-T)$.
\end{Lemma}

It follows from curvature condition of Theorem \ref{main result}, equality \eqref{cal of curvature 0} and Lemma \ref{lem:main} that we have
\begin{align}
&\eta \sqrt{-1}\Theta_{\tilde{h}}-\sqrt{-1}\partial \bar{\partial}
\eta\otimes\text{Id}_E-\sqrt{-1}g\partial\eta \wedge\bar{\partial}\eta\otimes\text{Id}_E\notag\\
=&S_k(-v_{t_{0},\varepsilon}(t))\sqrt{-1}\Theta_{h}
+S_k(-v_{t_{0},\varepsilon}(t))\sqrt{-1}\partial\bar{\partial}\psi\otimes\text{Id}_E\notag\\
+&\big(v'_{t_0,\epsilon}(\psi)\sqrt{-1}\partial\bar{\partial}(\psi)+
v''_{t_0,\epsilon}(\psi)\sqrt{-1}\partial(\psi)\wedge\bar{\partial}(\psi)\big)\otimes\text{Id}_E\notag\\
\ge &S_k\big(-\psi\big) v'_{t_0,\epsilon}(\psi)\bigg(\sqrt{-1}\Theta_{h}+\sqrt{-1}\partial\bar{\partial}\psi\bigg)\notag\\
+&v'_{t_0,\epsilon}(\psi)S_k(-\psi)\frac{1}{S_k(-\psi)}\sqrt{-1}\partial\bar{\partial}\psi\otimes \text{Id}_E+v''_{t_0,\epsilon}(\psi)\sqrt{-1}\big(\partial\psi\wedge\bar{\partial}\psi\big)\otimes\text{Id}_E\notag\\
=&S_k\big(-\psi\big) v'_{t_0,\epsilon}(\psi)\bigg(\sqrt{-1}\Theta_{h}+\sqrt{-1}\partial\bar{\partial}\psi+\frac{1}{S_k(-\psi)}\sqrt{-1}\partial\bar{\partial}\psi\bigg)\notag\\
+&v''_{t_0,\epsilon}(\psi)\sqrt{-1}\big(\partial\psi\wedge\bar{\partial}\psi\big)\otimes\text{Id}_E\notag\\
\ge &v''_{t_0,\epsilon}(\psi)\sqrt{-1}\big(\partial\psi\wedge\bar{\partial}\psi\big)\otimes\text{Id}_E.
\label{cal of curvature}
\end{align}

Then by \eqref{cal of curvature}, we have
\begin{equation}\label{step3 curvature}
\begin{split}
B
\ge
v''_{t_0,\epsilon}(\psi)[\sqrt{-1}\big(\partial\psi\wedge\bar{\partial}\psi\big)\otimes\text{Id}_E, \Lambda_{\omega}]
\end{split}
\end{equation}
holds on $M_k\backslash \sum$.

Let $\lambda_{t_0}:=D''[(1-v'_{t_0,\epsilon}(\psi))\tilde{f}]$. Then we know that $\lambda_{t_0}$ is well defined on $M_k$, $D''\lambda_{t_0}=0$ and
\begin{equation}\nonumber
\begin{split}
\lambda_{t_0}&=-v''_{t_0,\epsilon}(\psi)\bar{\partial}\psi\wedge\tilde{f}
+\big(1-v'_{t_0,\epsilon}(\psi)\big)D''\tilde{f}\\
&=\lambda_{1,t_0}+\lambda_{2,t_0},
\end{split}
\end{equation}
where $\lambda_{1,t_0}:=-v''_{t_0,\epsilon}(\psi)\bar{\partial}\psi\wedge\tilde{f}$ and
$\lambda_{2,t_0}:=\big(1-v'_{t_0,\epsilon}(\psi)\big)D''\tilde{f}$. Note that
$$\text{supp}\lambda_{1,t_0}\subset \{-t_0-1+\epsilon<\psi<-t_0-\epsilon\}$$
and
$$\text{supp}\lambda_{2,t_0}\subset \{\psi<-t_0-\epsilon\}.$$
It follows from inequality \eqref{step3 curvature} that we have
\begin{equation}\nonumber
\begin{split}
&\langle B^{-1}\lambda_{1,t_0},\lambda_{1,t_0} \rangle_{\omega,\tilde{h}}|_{M_k\backslash \sum}\\
\le &v''_{t_0,\epsilon}(\psi)|\tilde{f}|^2_{\omega,h}e^{-\psi-\phi}.
\end{split}
\end{equation}
Then we know that
\begin{equation}\label{step3 I1 1}
\begin{split}
&\int_{M_k\backslash\sum}\langle B^{-1}\lambda_{1,t_0},\lambda_{1,t_0} \rangle_{\tilde{h}}dV_{M,\omega}\\
\le &\int_{M_k\backslash\sum}v''_{t_0,\epsilon}(\psi)|\tilde{f}|^2_{\omega,h}e^{-\psi-\phi}dV_{M,\omega}\\
\le& I_{1,t_0,\epsilon}:=\sup_{M_k} e^{-\phi}\int_{M_k}v''_{t_0,\epsilon}(\psi)|\tilde{f}|^2_{\omega,h}e^{-\psi}dV_{M,\omega}.
\end{split}
\end{equation}

Denote $$I_{1,t_0}:=(\sup_{t\ge t_0} e^{-u_k(t)}) \int_{M_k\cap\{-t_0-1<\psi<-t_0\}}|\tilde{f}|^2_{\omega,h}e^{-\psi}dV_{M,\omega}.$$

It follows from  $\tilde{f}$ is a smooth extension of $f$ from $Y_0$ to $M$ and the definition of $dV_M[\psi]$ and $u(t)$ that we know
\begin{equation}\label{step3 I1 2}
\begin{split}
&\limsup_{t_0\to +\infty} I_{1,t_0}\\
\le& \big(\frac{1}{\delta}c_k(T)e^{-T}+\int^{+\infty}_T c_k(t_1)e^{-t_1}dt_1\big) \limsup_{t_0\to +\infty} \int_{M_k\cap\{-t_0-1<\psi<-t_0\}}|\tilde{f}|^2_{\omega,h}e^{-\psi}dV_{M,\omega}\\
\le & \big(\frac{1}{\delta}c_k(T)e^{-T}+\int^{+\infty}_T c_k(t_1)e^{-t_1}dt_1\big)
\int_{Y^0\cap M_k}|f|^2_{\omega,h}dV_{M,\omega}[\psi]\\
\le & \big(\frac{1}{\delta}c_k(T)e^{-T}+\int^{+\infty}_T c_k(t_1)e^{-t_1}dt_1\big)
\int_{Y^0}|f|^2_{\omega,h}dV_{M,\omega}[\psi].
\end{split}
\end{equation}

By inequalities \eqref{step3 I1 1} and \eqref{step3 I1 2}, we have
\begin{equation}\label{step3 I1m't0 2}
\begin{split}
&\limsup_{t_0\to +\infty}  I_{1,t_0,\epsilon}\\
\le &\frac{1}{1-4\epsilon}\limsup_{t_0\to +\infty}  I_{1,t_0}\\
\le  &\frac{1}{1-4\epsilon}\big(\frac{1}{\delta}c_k(T)e^{-T}+\int^{+\infty}_T c_k(t_1)e^{-t_1}dt_1\big)
\int_{Y^0}|f|^2_{\omega,h}dV_{M,\omega}[\psi].
\end{split}
\end{equation}

Denote
\begin{equation}\label{step3 I2 1}
\begin{split}
I_{2,t_0}:=&\int_{M_k\backslash \sum}\langle\lambda_{2,t_0},\lambda_{2,t_0} \rangle_{\tilde{h}}dV_{M,\omega}\\
\le&\int_{M_k\cap\{\psi<-t_0-\epsilon\}}|D'' \tilde{f}_{t_0}|^2_{\omega,h}e^{-\psi-\phi}dV_{M,\omega}\\
\le & \big(\sup_{t\ge t_0} e^{-u(t)}\big)
\int_{M_k\cap\{\psi<-t_0\}}|D'' \tilde{f}_{t_0}|^2_{\omega,h}e^{-\psi}dV_{M,\omega}.
\end{split}
\end{equation}

It follows from equality \eqref{d-baroftildef} and Cauchy-Schwarz inequality that when $t_0$ is big enough,
\begin{equation}\label{step3 I2 2}
\begin{split}
I_{2,t_0}\le& C_8\sum_{1\le i,j\le N}
\int_{U_i\cap U_j\cap\{\psi<-t_0\}}|\tilde{f}_{i}-\tilde{f}_{j}|^2_{\omega,h}e^{-\psi}dV_{M,\omega},
\end{split}
\end{equation}
where $C_8>0$ is a real number independent of $t_0$.

For any $1\le i,j\le N$, we denote
\begin{equation}\nonumber
\begin{split}
I_{i,j,t_0}:=\int_{U_i\cap U_j\cap\{\psi<-t_0\}}|\tilde{f}_{i}-\tilde{f}_{j}|^2_{\omega,h}e^{-\psi}dV_{M,\omega}.
\end{split}
\end{equation}
It follows from equality \eqref{fi-fj} and dominated convergence theorem that we have
\begin{equation}\nonumber
\begin{split}
\lim_{t_0\to +\infty}I_{i,j,t_0}=\lim_{t_0\to +\infty}\int_{U_i\cap U_j\cap\{\psi<-t_0\}}|\tilde{f}_{i}-\tilde{f}_{j}|^2_{\omega,h}e^{-\psi}dV_{M,\omega}=0.
\end{split}
\end{equation}
Hence we have
\begin{equation}\label{step3 I2 3}
\begin{split}
\lim_{t_0\to +\infty}I_{2,t_0}=0.
\end{split}
\end{equation}

\

\textbf{Step 5: solving $\bar{\partial}-$equation with error term.}

\

Given $\tau>0$,
note that $$\langle a_1+a_2,a_1+a_2\rangle\le (1+\tau)\langle a_1,a_1\rangle+(1+\frac{1}{\tau})\langle a_2,a_2\rangle$$
holds for any $a_1,a_2$ in an inner product space $(H,\langle \cdot,\cdot \rangle)$.
It follows from inequality \eqref{step3 curvature} that on $M_k\backslash \sum$, for any $\tau>0$, we have
\begin{equation}\label{Step 4 1}
\begin{split}
&\int_{M_k\backslash \sum} \langle \big(B+\sqrt{I_{2,t_0}}\text{Id}_E\big)^{-1}\lambda_{t_0},\lambda_{t_0} \rangle_{\omega,\tilde{h}}dV_{M,\omega}\\
\le & \int_{M_k\backslash\sum}(1+\tau) \langle \big(B+\sqrt{I_{2,t_0}}\text{Id}_E\big)^{-1}\lambda_{1,t_0},\lambda_{1,t_0} \rangle_{\omega,\tilde{h}}dV_{M,\omega}\\
&+\int_{M_k\backslash \sum}(1+\frac{1}{\tau}) \langle \big(B+\sqrt{I_{2,t_0}}\text{Id}_E\big)^{-1}\lambda_{2,t_0},\lambda_{2,t_0} \rangle_{\omega,\tilde{h}}dV_{M,\omega}\\
\le &(1+\tau)\int_{M_k\backslash\sum} \langle B^{-1}\lambda_{1,t_0},\lambda_{1,t_0} \rangle_{\omega,\tilde{h}}dV_{M,\omega}\\
+&(1+\frac{1}{\tau})\int_{M_k\backslash\sum} \langle \frac{1}{\sqrt{I_{2,t_0}}} \lambda_{2,t_0},\lambda_{2,t_0} \rangle_{\omega,\tilde{h}}dV_{M,\omega}\\
= & (1+\tau) I_{1,t_0,\epsilon}+(1+\frac{1}{\tau})\frac{1}{\sqrt{I_{2,t_0}}}I_{2,t_0}\\
\le &(1+\tau) I_{1,t_0,\epsilon}+(1+\frac{1}{\tau})\sqrt{I_{2,t_0}},
\end{split}
\end{equation}
By inequalities \eqref{step3 I1m't0 2} and \eqref{step3 I2 3}, we know that  $(1+\tau) I_{1,t_0,\epsilon}+(1+\frac{1}{\tau})\sqrt{I_{2,t_0}}$ is finite.

From now on, we fix some $\epsilon \in (0,\frac{1}{8})$.
 Then by Lemma \ref{d-bar equation with error term}, there exist $u_{k,t_0,\epsilon}\in L^2(M_k\backslash \sum, K_M\otimes E, \omega\otimes\tilde{h})$ and $\eta_{k,t_0,\epsilon}\in L^2(M_k\backslash \sum, \wedge^{n,1}T^*M\otimes E, \omega\otimes\tilde{h})$ such that
\begin{equation}\label{step 4 d-bar equation}
D''u_{k,t_0,\epsilon}+(I_{2,t_0})^{\frac{1}{4}}\eta_{k,t_0,\epsilon}=\lambda_{t_0}
\end{equation}
and
\begin{equation}\label{step 4 estimate}
\begin{split}
&\int_{M_k\backslash \sum}(\eta+g^{-1})^{-1}|u_{k,t_0,\epsilon}|^2_{\omega,\tilde{h}}dV_{M,\omega}
+\int_{M_k\backslash \sum}|\eta_{k,t_0,\epsilon}|^2_{\omega,\tilde{h}}dV_{M,\omega}\\
\le & (1+\tau) I_{1,t_0,\epsilon}+(1+\frac{1}{\tau})\sqrt{I_{2,t_0}}<+\infty.
\end{split}
\end{equation}
By definition, $(\eta+g^{-1})^{-1}=c_k(-v_{t_0,\epsilon}(\psi))e^{v_{t_0,\epsilon}(\psi)}e^{\phi}$. It follows from inequality \eqref{step 4 estimate} that
\begin{equation}\label{step 4 estimate 1}
\begin{split}
&\int_{M_k\backslash \sum}c_k\big(-v_{t_0,\epsilon}(\psi)\big)e^{v_{t_0,\epsilon}(\psi)-\psi}|u_{k,t_0,\epsilon}|^2_{\omega,h}dV_{M,\omega}\\
\le & (1+\tau) I_{1,t_0,\epsilon}+(1+\frac{1}{\tau})\sqrt{I_{2,t_0}},
\end{split}
\end{equation}
and
\begin{equation}\label{step 4 estimate 2}
\begin{split}
&\int_{M_k\backslash \sum}|\eta_{k,t_0,\epsilon}|^2_{\omega,h}e^{-\psi-\phi}dV_{M,\omega}\\
\le  &(1+\tau) I_{1,t_0,\epsilon}+(1+\frac{1}{\tau})\sqrt{I_{2,t_0}}.
\end{split}
\end{equation}

Note that $v_{t_0,\epsilon}(\psi)$ is bounded on $M_k$ and $c_k(t)e^{-t}$ is decreasing near $+\infty$.  We know that $c(-v_{t_0,\epsilon}(\psi))e^{v_{t_0,\epsilon}(\psi)}$ has positive lower bound on $M_k$. We also have $e^{-\phi}=e^{-u(-v_{t_0,\epsilon}(\psi))}$ has positive lower bound on $M_k$. As $\psi$ is upper-bounded on $M_k$, $e^{-\psi}$ also have positive lower bound on $M_k$. By inequalities \eqref{step 4 estimate 1} and \eqref{step 4 estimate 2}, we know that
$$u_{k,t_0,\epsilon}\in L^2(M_k\backslash \sum, K_M\otimes E, \omega\otimes h)$$
and $$\eta_{k,t_0,\epsilon}\in L^2(M_k\backslash \sum, \wedge^{n,1}T^*M\otimes E, \omega\otimes h).$$
By Lemma \ref{extension of equality} and equality \eqref{step 4 d-bar equation}, we know that
\begin{equation}\label{step 4 d-bar equation 2}
D''u_{k,t_0,\epsilon}+(I_{2,t_0})^{\frac{1}{4}}\eta_{k,t_0,\epsilon}=\lambda_{t_0}
\end{equation}
holds on $M_k$.

Recall that $\lambda_{t_0}:=D''[(1-v'_{t_0,\epsilon}(\psi))\tilde{f}]$. Denote $F_{k,t_0,\epsilon}:=\lambda_{t_0}-u_{k,t_0,\epsilon}$. Then we have
\begin{equation}\label{step 4 d-bar equation 3}
D''F_{k,t_0,\epsilon}=(I_{2,t_0})^{\frac{1}{4}}\eta_{k,t_0,\epsilon}
\end{equation}
holds on $M_k$.
It follows from $\sum$ is a set of measure zero and inequalities \eqref{step 4 estimate 1} and \eqref{step 4 estimate 2} that we have
\begin{equation}\label{step 4 estimate 3}
\begin{split}
&\int_{M_k}c_k\big(-v_{t_0,\epsilon}(\psi)\big)e^{v_{t_0,\epsilon}(\psi)-\psi}
|u_{k,t_0,\epsilon}|^2_{\omega,h}dV_{M,\omega}\\
=&\int_{M_k}c_k\big(-v_{t_0,\epsilon}(\psi)\big)e^{v_{t_0,\epsilon}(\psi)-\psi}
|F_{k,t_0,\epsilon}-\big((1-v'_{t_0,\epsilon}(\psi))\tilde{f}\big)|^2_{\omega,h}dV_{M,\omega}\\
\le & (1+\tau) I_{1,t_0,\epsilon}+(1+\frac{1}{\tau})\sqrt{I_{2,t_0}},
\end{split}
\end{equation}
and
\begin{equation}\label{step 4 estimate 4}
\begin{split}
&\int_{M_k}|\eta_{k,t_0,\epsilon}|^2_{\omega,h}e^{-\psi-\phi}dV_{M,\omega}\\
\le  &(1+\tau) I_{1,t_0,\epsilon}+(1+\frac{1}{\tau})\sqrt{I_{2,t_0}}.
\end{split}
\end{equation}

\

\textbf{Step 6: when $t_0\to +\infty$.}

\

We note that $k$ is fixed in this step. It follows from inequality \eqref{step 4 estimate 3}, $c_k(t)e^{-t}$ is decreasing with respect to $t$ near $+\infty$ and $v_{t_0,\epsilon}(\psi)\ge \psi$ that, when $t_0$ is big enough, we have

\begin{equation}\label{step 7 estimate 1}
\begin{split}
&\int_{M_k}c_k(-\psi)
|F_{k,t_0,\epsilon}-\big((1-v'_{t_0,\epsilon}(\psi))\tilde{f}\big)|^2_{\omega,h}dV_{M,\omega}\\
\le & (1+\tau) I_{1,t_0,\epsilon}+(1+\frac{1}{\tau})\sqrt{I_{2,t_0}}.
\end{split}
\end{equation}
Then we have
\begin{equation}\label{step 7 estimate 2}
\begin{split}
&\int_{M_k}c_k(-\psi)
|F_{k,t_0,\epsilon}|^2_{\omega,h}dV_{M,\omega}\\
\le &(1+\tau)\int_{M_k}c_k(-\psi)
|F_{k,t_0,\epsilon}-\big((1-v'_{t_0,\epsilon}(\psi))\tilde{f}\big)|^2_{\omega,h}dV_{M,\omega}\\
+&(1+\frac{1}{\tau})\int_{M_k}c_k(-\psi)
|\big((1-v'_{t_0,\epsilon}(\psi))\tilde{f}\big)|^2_{\omega,h}dV_{M,\omega}\\
\le & (1+\tau)^2 I_{1,t_0,\epsilon}+(1+\tau)(1+\frac{1}{\tau})\sqrt{I_{2,t_0}}\\
+&(1+\frac{1}{\tau})C_{t_0},
\end{split}
\end{equation}
where $C_{t_0}:=\int_{M_k}c_k(-\psi)
|\big((1-v'_{t_0,\epsilon}(\psi))\tilde{f}\big)|^2_{\omega,h}dV_{M,\omega}$. Note that $(1-v'_{t_0,\epsilon}(\psi))\le \mathbb{I}_{\{\psi<-t_0\}}$. We have
\begin{equation}\label{cal of ct0}
\begin{split}
C_{t_0}=&\int_{M_k}c_k(-\psi)
|\big((1-v'_{t_0,\epsilon}(\psi))\tilde{f}\big)|^2_{\omega,h}dV_{M,\omega}\\
\le & \int_{M_k}c_k(-\psi)\mathbb{I}_{\{\psi<-t_0\}}
|\tilde{f}|^2_{\omega,h}dV_{M,\omega}\\
= & \sum_{j=0}^{+\infty}\int_{M_k}c_k(-\psi)\mathbb{I}_{\{-t_0-1-j\le\psi<-t_0-j\}}
|\tilde{f}|^2_{\omega,h}dV_{M,\omega}.
\end{split}
\end{equation}
As $c_k(t)e^{-t}$ is decreasing with respect to $t$ near $+\infty$, we take $t_0$ big enough, such that $c_k(t)e^{-t}$ is decreasing with respect to $t$ on $[t_0-1,+\infty)$. Then on $\{t_0+j<-\psi\le t_0+1+j\}$, we have $c_k(-\psi)e^{\psi}\le c(t_0+j)e^{-t_0-j}$. Combining with inequality \eqref{cal of ct0}, we have
\begin{equation}\label{cal of ct0 1}
\begin{split}
C_{t_0}\le & \sum_{j=0}^{+\infty}\int_{M_k}c_k(-\psi)\mathbb{I}_{\{-t_0-1-j\le\psi<-t_0-j\}}
|\tilde{f}|^2_{\omega,h}dV_{M,\omega}\\
\le & \sum_{j=0}^{+\infty}\int_{M_k}c_k(t_0+j)e^{-t_0-j}e^{-\psi}\mathbb{I}_{\{-t_0-1-j\le\psi<-t_0-j\}}
|\tilde{f}|^2_{\omega,h}dV_{M,\omega}\\
=& \sum_{j=0}^{+\infty}c_k(t_0+j)e^{-t_0-j}\int_{M_k}\mathbb{I}_{\{-t_0-1-j\le\psi<-t_0-j\}}
|\tilde{f}|^2_{\omega,h}e^{-\psi}dV_{M,\omega}.
\end{split}
\end{equation}
It follows from condition \eqref{mainth:ohsawa measure finite} that
$$\limsup_{t\rightarrow+\infty}\int_{\{-t-1<\psi<-t\}}| \tilde{f}|^2_{\omega,h}e^{-\psi}dV_{M,\omega}<+\infty.$$
Hence, when $t_0$ is big enough, we have
\begin{equation}\nonumber
\begin{split}
\sup_{j\ge 0}\int_{M_k}\mathbb{I}_{\{-t_0-1-j\le\psi<-t_0-j\}}
|\tilde{f}|^2_{\omega,h}e^{-\psi}dV_{M,\omega}\le C_1,
\end{split}
\end{equation}
for some constant $C_1>0$ independent of $j$. Then we know
\begin{equation}\nonumber
\begin{split}
C_{t_0}\le
& C_1\sum_{j=0}^{+\infty}c_k(t_0+j)e^{-t_0-j}\le C_1\int_{t_0-1}^{+\infty}c_k(t)e^{-t}dt,
\end{split}
\end{equation}
and  we have
\begin{equation}\label{cal of ct0 3}
\begin{split}
\lim_{t_0\to +\infty}C_{t_0}=0.
\end{split}
\end{equation}

By inequalities \eqref{step 7 estimate 2}, \eqref{step3 I1m't0 2}, \eqref{step3 I2 3} and \eqref{cal of ct0 3}, when $t_0$ is big enough, we have,
$$\sup_{t_0}\int_{M_k}c_k(-\psi)
|F_{k,t_0,\epsilon}|^2_{\omega,h}dV_{M,\omega}<+\infty.$$
Since the closed unit ball of Hilbert space is
weakly compact, we know that there exists a subsequence of $\{F_{k,t_0,\epsilon}\}_{t_0}$ (also denoted by $\{F_{k,t_0,\epsilon}\}$) weakly
convergent to some $F_{k,\epsilon}$ in $L^2(M_k, K_M\otimes E, c_k(-\psi)\omega\otimes h)$ as $t_0\to+\infty$. Again by inequalities \eqref{step 7 estimate 2}, \eqref{step3 I1m't0 2}, \eqref{step3 I2 3} and \eqref{cal of ct0 3}, we have
\begin{equation}\nonumber
\begin{split}
&\int_{M_k}c_k(-\psi)
|F_{k,\epsilon}|^2_{\omega,h}dV_{M,\omega}\\
\le &\liminf_{t_0\to +\infty}\int_{M_k}c_k(-\psi)
|F_{k,t_0,\epsilon}|^2_{\omega,h}dV_{M,\omega}\\
\le & \limsup_{t_0\to +\infty}\bigg((1+\tau)^2 I_{1,t_0,\epsilon}+(1+\tau)(1+\frac{1}{\tau})\sqrt{I_{2,t_0}}+(1+\frac{1}{\tau})C_{t_0}\bigg)\\
=&\frac{(1+\tau)^2}{1-4\epsilon}\big(\frac{1}{\delta}c_k(T)e^{-T}+\int^{+\infty}_T c_k(t_1)e^{-t_1}dt_1\big)
\int_{Y^0}|f|^2_{\omega,h}dV_{M,\omega}[\psi].
\end{split}
\end{equation}

By the construction of $e^{-\phi}$ and $\sup_{M_k}\psi<-T-8\epsilon_k$, we know that $e^{-\phi}=e^{-u_k(-v_{t_0,\epsilon}(\psi))}$ has a uniformly  positive lower bounded with respect to $t_0$. It follows from  $\psi$ is upper bounded on $M_k$, inequalities \eqref{step 4 estimate 4}, \eqref{step3 I1m't0 2} and \eqref{step3 I2 3} that we have
$$\sup_{t_0}\int_{M_k}|\eta_{k,t_0,\epsilon}|^2_{\omega,h}dV_{M,\omega}<+\infty.$$
Since the closed unit ball of Hilbert space is
weakly compact, we can extract a subsequence of $\{\eta_{k,t_0,\epsilon}\}_{t_0}$ (also denoted by $\{v_{k,t_0,\epsilon}\}_{t_0}$) weakly convergent to $v_{k,\epsilon}$ in
$L^2(M_k, \wedge^{n,1}T^*M\otimes E, \omega\otimes h)$ as $t_0\to+\infty$. As $\lim_{t_0\to +\infty}I_{2,t_0}=0$, we know $I_{2,t_0}$ is uniformly upper bounded with respect to $t_0$. It follows from Lemma \ref{weakly convergence} that
$(I_{2,t_0})^{\frac{1}{4}}\eta_{k,t_0,\epsilon}$ is weakly convergent to $0$ in
$L^2(M_k, \wedge^{n,1}T^*M\otimes E, \omega\otimes h)$ as $t_0\to+\infty$. Hence $(I_{2,t_0})^{\frac{1}{4}}\eta_{k,t_0,\epsilon}$ is weakly convergent to $0$ in
$L^2_{\text{loc}}(M_k, \wedge^{n,1}T^*M\otimes E, \omega\otimes h)$ as $t_0\to+\infty$

It follows from $\psi$ is smooth on $M_k\backslash \sum$, $c_k(t)$ is smooth function on $[T+4\epsilon_k,+\infty)$ and $\{F_{k,t_0,\epsilon}\}$  weakly converges to
$\{F_{k,\epsilon}\}$ in $L^2(M_k, K_M\otimes E, c_k(-\psi) \omega\otimes h)$ as $t_0\to +\infty$ that we have $\{F_{k,t_0,\epsilon}\}$ also weakly converges to
$\{F_{k,\epsilon}\}$ in $L^2_{\text{loc}}(M_k\backslash \sum, K_M\otimes E, \omega\otimes h)$ as $t_0\to +\infty$.

Then it follows from equality \eqref{step 4 d-bar equation 3} that, when $t_0\to +\infty$, we have
\begin{equation}\label{step 7 d-bar equation}
D''F_{k,\epsilon}=0 \text{\ holds on } M_k\backslash \sum.
\end{equation}
 Hence $F_{k,\epsilon}$ is an $E$-valued holomorphic $(n,0)$-form on $M_k\backslash \sum$ which satisfies
\begin{equation}\label{step 7 estimate 3}
\begin{split}
&\int_{M_k}c_k(-\psi)
|F_{k,\epsilon}|^2_{\omega,h}dV_{M,\omega}\\
\le &\frac{(1+\tau)^2}{1-4\epsilon}\big(\frac{1}{\delta}c_k(T)e^{-T}+\int^{+\infty}_T c_k(t_1)e^{-t_1}dt_1\big)
\int_{Y^0}|f|^2_{\omega,h}dV_{M,\omega}[\psi].
\end{split}
\end{equation}
\

\textbf{Step 7: solving $\bar{\partial}$-equation locally.}

\

In this step, we prove that $F_{k,\epsilon}$ is actually  a holomorphic extension of $f$ from $Y^0\cap M_k$ to $M_k$.

Let $x\in M_k\cap Y^0$ be any point.  Let $\tilde{U}_x\Subset  M_k$ be a local coordinate ball which is centered at $x$. We assume that $E|_{\tilde{U}_x}$ is trivial vector bundle.

Note that  we have
$$D''u_{k,t_0,\epsilon}+(I_{2,t_0})^{\frac{1}{4}}\eta_{k,t_0,\epsilon}=
D''[\big(1-v'_{t_0,\epsilon}(\psi)\big)\tilde{f}].$$
It follows from inequality \eqref{step 4 estimate 4} and the construction of $e^{-\phi}$ that we have
\begin{equation}\label{Step 8 formula 1}
\int_{\tilde{U}_x}|\eta_{k,t_0,\epsilon}|^2_{h}e^{-\psi}\le C_{2},
\end{equation}
where  $C_{2}>0$ is a positive number independent of $t_0$.

Note that $D''\big((I_{2,t_0})^{\frac{1}{4}}\eta_{k,t_0,\epsilon}\big)=0$. It follows from Lemma \ref{hormander} that there exists an $E$-valued $(n,0)$-form $s_{k,t_0,\epsilon}\in L^2(\tilde{U}_x, K_M\otimes E,he^{-\psi})$ such that $D''s_{k,t_0,\epsilon}=(I_{2,t_0})^{\frac{1}{4}}\eta_{k,t_0,\epsilon}$ and
\begin{equation}\label{Step 8 formula 2}
\int_{\tilde{U}_x}|s_{k,t_0,\epsilon}|^2_{h}e^{-\psi}\le C_{3}\int_{\tilde{U}_x}|(I_{2,t_0})^{\frac{1}{4}}\eta_{k,t_0,\epsilon}|^2_{h}e^{-\psi}
\le C_{3}C_{2}\sqrt{I_{2,t_0}},
\end{equation}
where $C_{3}>0$ is a positive number independent of $t_0$. Hence we have
\begin{equation}\label{Step 8 formula 3}
\int_{\tilde{U}_x}|s_{k,t_0,\epsilon}|^2_{h}\le C_{4}\sqrt{I_{2,t_0}},
\end{equation}
where $C_{4}>0$ is a positive number independent of $t_0$.

Now define $G_{k,t_0,\epsilon}:=-u_{k,t_0,\epsilon}-s_{k,t_0,\epsilon}
+\big(1-v'_{t_0,\epsilon}(\psi)\big)\tilde{f}$ on $\tilde{U}_x$. Then we know that
$G_{k,t_0,\epsilon}=F_{k,t_0,\epsilon}-s_{k,t_0,\epsilon}$ and $D''G_{k,t_0,\epsilon}=0$. Hence $G_{k,t_0,\epsilon}$ is holomorphic on $\tilde{U}_x$ and we know that $u_{k,t_0,\epsilon}+s_{k,t_0,\epsilon}$ is smooth on $\tilde{U}_x$ .

It follows from $c_k(t)e^{-t}$ is decreasing with respect to $t$ near $+\infty$ and inequalities \eqref{step 7 estimate 2} and \eqref{Step 8 formula 2} that we have
\begin{equation}\label{Step 8 formula 4}
\begin{split}
&\int_{\tilde{U}_x}c_k(-\psi)|G_{k,t_0,\epsilon}|^2_{h}\\
\le &2\tilde{C_{5}}\int_{\tilde{U}_x}c_k(-\psi)|F_{k,t_0,\epsilon}|^2_{h}
+2\tilde{C_{5}}\int_{\tilde{U}_x}|s_{k,t_0,\epsilon}|^2_{h}e^{-\psi}\le C_{5},
\end{split}
\end{equation}
where $\tilde{C_{5}},C_{5}>0$ are positive numbers independent of $t_0$.

It follows from inequality \eqref{step 4 estimate 3}, the construction of $v_{t_0,\epsilon}(\psi)$ and $c_k(t)e^{-t}$ is decreasing with respect to $t$ near $+\infty$ that we have
\begin{equation}\label{Step 8 formula 5}
\int_{\tilde{U}_x}|u_{k,t_0,\epsilon}|^2_{h}e^{-\psi}
\le C_{6,t_0},
\end{equation}
where $C_{6,t_0}>0$ is a sequence of positive number depending on $t_0$. Then, by inequalities \eqref{Step 8 formula 2} and \eqref{Step 8 formula 5}, we have
\begin{equation}\label{Step 8 formula 6}
\int_{\tilde{U}_x}|u_{k,t_0,\epsilon}+s_{k,t_0,\epsilon}|^2_{\tilde{h}}e^{-\psi}
\le 2C_{6,t_0}+2C_{3}C_{2}\sqrt{I_{2,t_0}}.
\end{equation}
Note that $e^{-\psi}$ is not integrable along $Y$ and $u_{k,t_0,\epsilon}+s_{k,t_0,\epsilon}$ is smooth on $\tilde{U}_x$. By \eqref{Step 8 formula 6}, we know that $u_{k,t_0,\epsilon}+s_{k,t_0,\epsilon}=0$ on $\tilde{U}_x\cap Y$ for any $t_0$. Hence $G_{k,t_0,\epsilon}
=\tilde{f}_{t_0}=f$ on $\tilde{U}_x\cap Y_0$ for any $t_0$.

It follows from inequality \eqref{Step 8 formula 3} that there exists a subsequence of $\{s_{k,t_0,\epsilon}\}$ \big(also denoted by $\{s_{k,t_0,\epsilon}\}$\big) weakly convergent to $0$ in  $L^2(\tilde{U}_x, K_M\otimes E,h)$ as $t_0\to +\infty$. Note that $\{F_{k,t_0,\epsilon}\}$ weakly converges to $F_{k,\epsilon}$ in  $L^2_{\text{loc}}(\tilde{U}_x\backslash \sum, K_M\otimes E,h)$ as $t_0\to +\infty$.  Hence we know that $\{G_{k,t_0,\epsilon}\}$ weakly converges to $F_{k,\epsilon}$ in  $L^2_{\text{loc}}(\tilde{U}_x\backslash \sum, K_M\otimes E,h)$as $t_0\to +\infty$.

It follows from inequality \eqref{Step 8 formula 4}  and Lemma \ref{l:converge} that we know there exists a subsequence of $\{G_{k,t_0,\epsilon}\}$ \big(also denoted by $\{G_{k,t_0,\epsilon}\}$\big) compactly convergent to an $E$-valued holomorphic $(n,0)$-form $G_{k,\epsilon}$ on $\tilde{U}_x$ as $t_0\to +\infty$. As $G_{k,t_0,\epsilon}
=f$ on $\tilde{U}_x\cap Y_0$ for any $t_0$, we know that $G_{k,\epsilon}=f$ on  $\tilde{U}_x\cap Y_0$.

As $\{G_{k,t_0,\epsilon}\}$  compactly converges to $G_{k,\epsilon}$ on $\tilde{U}_x$ as $t_0\to +\infty$ and  $\{G_{k,t_0,\epsilon}\}$ weakly converges to $F_{k,\epsilon}$ in  $L^2_{\text{loc}}(\tilde{U}_x\backslash \sum, K_M\otimes E,h)$ as $t_0\to +\infty$, by the uniqueness of weak limit, we know that $G_{k,\epsilon}=F_{k,\epsilon}$ on any relatively compact open subset of $\tilde{U}_x$. Note that $G_{k,,\epsilon}$ is holomorphic on $\tilde{U}_x$ and $F_{k,\epsilon}$ is holomorphic on $\tilde{U}_x\backslash \sum$, we have $F_{k,\epsilon}\equiv G_{k,\epsilon}$ on $\tilde{U}_x\backslash \sum$, and we know that $F_{k,\epsilon}$ can extended to an $E$-valued holomorphic $(n,0)$-form on $\tilde{U}_x$ which equals to $G_{k,\epsilon}$. As $G_{k,\epsilon}=f$ on  $\tilde{U}_x\cap Y_0$, we know that $F_{k,\epsilon}=f$ on  $\tilde{U}_x\cap Y_0$. Since $x$ is arbitrarily chosen, we know that $F_{k,\epsilon}$ is  holomorphic on $M_k$ and $F_{k,\epsilon}=f$ on  $M_k\cap Y_0$.

\

\textbf{Step 8: end of the proof.}

\

Now we have a family of $E$-valued holomorphic $(n,0)$-forms $F_{k,\epsilon}$ on $M_k$ such that $F_{k,\epsilon}=f$ on $M_k\cap Y_0$ and
\begin{equation}\label{step 8 estimate 0}
\begin{split}
&\int_{M_k} c_k(-\psi)|F_{k,\epsilon}|^2_{\omega,h}dV_{M,\omega}\\
\le&\frac{(1+\tau)^2}{1-4\epsilon}\big(\frac{1}{\delta}c_k(T)e^{-T}+\int^{+\infty}_T c_k(t_1)e^{-t_1}dt_1\big)
\int_{Y^0}|f|^2_{\omega,h}dV_{M,\omega}[\psi].
\end{split}
\end{equation}
Recall that $\epsilon\in (0,\frac{1}{8})$. By inequality \eqref{step 8 estimate 0}, we have
\begin{equation}\label{step 11 estimate 0}
\sup_{\epsilon\in (0,\frac{1}{8})}\int_{M_k} c_k(-\psi)|F_{k,\epsilon}|^2_{\omega,h}dV_{M,\omega}<+\infty.
\end{equation}

For any compact subset $K\subset M_k\backslash \sum=\{\psi=-\infty\}$, as $\psi$ is smooth on $M_k\backslash \sum$, we know that $\psi$ is upper and lower bounded on $K$. As $c_k(t)$ is continuous on $[T,+\infty)$, we have $c_{k}(-\psi)$ is uniformly lower bounded on $K$.
 It follows from Lemma \ref{l:converge} and inequality \eqref{step 11 estimate 0} that there exists a subsequence of $\{F_{k,\epsilon}\}_{\epsilon}$ (also denoted by $\{F_{k,\epsilon}\}_{\epsilon}$) compactly convergent to an $E$-valued holomorphic $(n,0)$-form  $F_k$ on $M_k$ as $\epsilon \to 0$. It follows from Fatou's lemma (let $\epsilon \to 0$) and inequality \eqref{step 8 estimate 0} that we have
\begin{equation}\nonumber
\begin{split}
&\int_{M_k} c_k(-\psi)|F_{k}|^2_{\omega,h}dV_{M,\omega}\\
\le&\liminf_{\epsilon \to 0}\int_{M_k} c_k(-\psi)|F_{k,\epsilon}|^2_{\omega,h}dV_{M,\omega}\\
\le&\liminf_{\epsilon \to 0}\frac{(1+\tau)^2}{1-4\epsilon}\big(\frac{1}{\delta}c_k(T)e^{-T}+\int^{+\infty}_T c_k(t_1)e^{-t_1}dt_1\big)
\int_{Y^0}|f|^2_{\omega,h}dV_{M,\omega}[\psi]\\
\le &(1+\tau)^2\big(\frac{1}{\delta}c_k(T)e^{-T}+\int^{+\infty}_T c_k(t_1)e^{-t_1}dt_1\big)
\int_{Y^0}|f|^2_{\omega,h}dV_{M,\omega}[\psi].
\end{split}
\end{equation}
Hence there exists a family of $E$-valued holomorphic $(n,0)$-forms $F_{k}$ on $M_k$ such that $F_{k}=f$ on $M_k\cap Y_0$ and
\begin{equation}\nonumber
\begin{split}
&\int_{M_k} c_k(-\psi)|F_{k}|^2_{\omega,h}dV_{M,\omega}\\
\le&(1+\tau)^2\big(\frac{1}{\delta}c_k(T)e^{-T}+\int^{+\infty}_T c_k(t_1)e^{-t_1}dt_1\big)
\int_{Y^0}|f|^2_{\omega,h}dV_{M,\omega}[\psi].
\end{split}
\end{equation}

Since $\tau>0$ is arbitrarily chosen and $c_k(T)e^{-T}=c(T)e^{-T}$, we have

\begin{equation}\label{continuous 1}
\begin{split}
&\int_{M_k} c_k(-\psi)|F_k|^2_{\omega,h}dV_{M,\omega}\\
\le& \big(\frac{1}{\delta}c(T)e^{-T}+\int^{+\infty}_T c_k(t_1)e^{-t_1}dt_1\big)
\int_{Y^0}|f|^2_{\omega,h}dV_{M,\omega}[\psi].
\end{split}
\end{equation}

Let $k_1> k$ be big enough. It follows from inequality \eqref{continuous 1}, $M_k\Subset M_{k_1}$ and $\frac{1}{\delta}c(T)e^{-T}+\int_{T}^{+\infty}c_k(t)e^{-t}dt$ converges to $\frac{1}{\delta}c(T)e^{-T}+\int_{T}^{+\infty}c(t)e^{-t}dt<+\infty$ as $k\to +\infty$ that we have
\begin{equation}\label{continuous uni estimate}
\begin{split}
\sup_{k_1}\int_{M_k} c_{k_1}(-\psi)|F_{k_1}|^2_{\omega,h}dV_{M,\omega}<+\infty.
\end{split}
\end{equation}

For any compact subset $K\subset M_k\backslash \sum=\{\psi=-\infty\}$, as $\psi$ is smooth on $M_k\backslash \sum$, we know that $\psi$ is upper and lower bounded on $K$. It follows from $c_{k_1}(t)$ are uniformly convergent to $c(t)$ on any compact subset of $(T,+\infty)$ and $c(t)$ is a positive continuous function on $[T,+\infty)$ that we know
$c_{k_1}(-\psi)$ is uniformly lower bounded on $K$. By Lemma \ref{l:converge} and inequality \eqref{continuous uni estimate}, we know that there exists a subsequence of $\{F_{k_1}\}_{{k_1}\in\mathbb{Z}^+}$ (also denoted by $\{F_{k_1}\}_{{k_1}\in\mathbb{Z}^+}$) compactly convergent to an $E$-valued holomorphic $(n,0)$-form  $\tilde{F}_k$ on $M_k$. It follows from Fatou's lemma (let $k_1\to +\infty$) and inequality \eqref{continuous 1} that we have
\begin{equation}\label{continuous 2}
\begin{split}
&\int_{M_k} c(-\psi)|\tilde{F}_k|^2_{\omega,h}dV_{M,\omega}\\
\le& \big(\frac{1}{\delta}c(T)e^{-T}+\int^{+\infty}_T c(t_1)e^{-t_1}dt_1\big)
\int_{Y^0}|f|^2_{\omega,h}dV_{M,\omega}[\psi].
\end{split}
\end{equation}
As $\{F_{k_1}\}_{{k_1}\in\mathbb{Z}^+}$  compactly convergent to an $E$-valued holomorphic $(n,0)$-form  $\tilde{F}_k$ on $M$, we know that $\tilde{F}_k=f$ on $M_k\cap Y_0$.

Again for any compact subset $K\subset M\backslash Y$, as $\psi$ is smooth on $M\backslash Y$, we know that $\psi$ is upper and lower bounded on $K$. It follows  $c(t)$ is a positive continuous function on $[T,+\infty)$ that we know
$c(-\psi)$ is uniformly lower bounded on $K$. By Lemma \ref{l:converge} and inequality \eqref{continuous 2}, we know that there exists a subsequence of $\{\tilde{F}_k\}_{k\in\mathbb{Z}^+}$ (also denoted by $\{\tilde{F}_k\}_{k\in\mathbb{Z}^+}$) compactly convergent to an $E$-valued holomorphic $(n,0)$-form  $F$ on $M$. It follows from Fatou's lemma (let $k\to +\infty$) and inequality \eqref{continuous 2} that we have
\begin{equation}\nonumber
\begin{split}
&\int_{M_k} c(-\psi)|F|^2_{\omega,h}dV_{M,\omega}\\
\le& \big(\frac{1}{\delta}c(T)e^{-T}+\int^{+\infty}_T c(t_1)e^{-t_1}dt_1\big)
\int_{Y^0}|f|^2_{\omega,h}dV_{M,\omega}[\psi].
\end{split}
\end{equation}
Letting $k\to +\infty$, by monotone convergence theorem, we have
\begin{equation}\nonumber
\begin{split}
&\int_{M} c(-\psi)|F|^2_{\omega,h}dV_{M,\omega}\\
\le& \big(\frac{1}{\delta}c(T)e^{-T}+\int^{+\infty}_T c(t_1)e^{-t_1}dt_1\big)
\int_{Y^0}|f|^2_{\omega,h}dV_{M,\omega}[\psi].
\end{split}
\end{equation}
As $\{\tilde{F}_k\}_{k\in\mathbb{Z}^+}$ is compactly convergent to an $E$-valued holomorphic $(n,0)$-form  $F$ on $M$, we know that $F=f$ on $M\cap Y_0$.

Theorem \ref{main result} has been proved.
\end{proof}


\vspace{.1in} {\em Acknowledgements}. The authors would like to thank Professor Xiangyu Zhou for his encouragement.
The authors would like to thank Dr. Shijie Bao for checking the manuscript.
The first author and the second author were supported by National Key R\&D Program of China 2021YFA1003100.
The first author was supported by NSFC-11825101, NSFC-11522101 and NSFC-11431013. The second author was supported by China Postdoctoral Science Foundation 2022T150687.

\bibliographystyle{references}
\bibliography{xbib}

\end{document}